\documentclass{article}
\usepackage{amsmath, amssymb, amsthm, epsf, graphicx}

\newtheorem{theorem}{Theorem}
\newtheorem{lemma}{Lemma}

\newtheorem{definition}{Definition}
\newtheorem{corollary}{Corollary}

\newtheorem{question}{Question}

\def\bb #1{ {\mathbb #1} }
\def\cc #1{ {\mathcal #1} }

\begin{document}
\title{On the zeta function of divisors for projective varieties \\ with higher rank divisor class group}

\author{C. Douglas Haessig \\
    University of Rochester, New York \\
    Department of Mathematics\\
    Rochester, NY 14627 \\
    chaessig@math.rochester.edu}

\date{\today}
\maketitle

\begin{abstract}
Given a projective variety $X$ defined over a finite field, the
zeta function of divisors attempts to count all irreducible,
codimension one subvarieties of $X$, each measured by their
projective degree. When the dimension of $X$ is greater than
one, this is a purely $p$-adic function, convergent on the open
unit disk. Four conjectures are expected to hold, the first of
which is $p$-adic meromorphic continuation to all of $\bb
C_p$.

When the divisor class group (divisors modulo linear
equivalence) of $X$ has rank one, then all four conjectures are
known to be true. In this paper, we discuss the higher rank
case. In particular, we prove a $p$-adic meromorphic
continuation theorem which applies to a large class of
varieties. Examples of such varieties are projective
nonsingular surfaces defined over a finite field (whose
effective monoid is finitely generated) and all projective
toric varieties (smooth or singular).
\end{abstract}

\section{Introduction}

This paper continues the study of the zeta function of
divisors, a function first proposed by Professor Daqing Wan in
\cite{WanZetaOfCycles}. Given a projective variety $X$ defined
over a finite field with a fixed embedding into projective
space, the zeta function of divisors $Z_{div}(X,T)$ is a
function generated by all irreducible subvarieties $P$ of
codimension one, where each subvariety is measured by its
projective degree:
\begin{equation}
Z_{div}(X,T) := \prod_{P \text{ irred}} \frac{1}{1-T^{deg(P)}}.
\end{equation}
Notice that when $X$ is a curve, this function agrees with the
classical Hasse-Weil zeta function since divisors are now points.
Unlike the Weil zeta function, which is always a rational
function in $T$ by Dwork's theorem, when the dimension of $X$
is greater than one $Z_{div}(X,T)$ is never a rational function
\cite{WanZetaOfCycles}. In fact it has zero radius of
convergence over the complex numbers. However, over the
$p$-adic numbers the situation is much better. An immediate
consequence of the definition is that $Z_{div}(X,T)$ is
$p$-adic analytic on the open unit disk. The study of whether
this function can be $p$-adic meromorphically continued makes
up the content of this paper.

Denote by $A^1(X)$ the divisor class group of $X$ (divisors
modulo linear equivalence) and let $A^+(X)$ be the image in
$A^1(X)$ of all effective divisors. We call $A^+(X)$ the
effective monoid. \emph{Assume the effective monoid is finitely
generated.} Then we have the following four conjectures.

\begin{description}
    \item[Conjecture I]($p$-adic Meromorphic continuation)
        \\ $Z_{div}(X,T)$ is $p$-adic meromorphic everywhere in $T$.
    \item[Conjecture II] (Order and rank)
        \\ $Z_{div}(X,T)$ has a pole at $T=1$ of order equal to the rank of
        $A^1(X)$.
    \item[Conjecture III] (Simplicity of zeros and poles)
        \\ Except for possibly finitely many, all zeros and poles are
        simple.
    \item[Conjecture IV] ($p$-adic Riemann hypothesis)
        \\ Adjoining the zeros and poles of $Z_{div}(X,T)$ to the
        field $\bb Q_p$ produces a finite extension field
        of $\bb Q_p$.
\end{description}
Conjectures I, II, and III were first stated in
\cite{WanZetaOfCycles}, while IV was stated in
\cite{HaessigWanDivisors}.

Previous studies of the zeta function of divisors have focused
on the case when $A^1(X)$ is of rank one. In fact, under this
condition the conjectures have been completely answered:

\begin{theorem}[\cite{WanZetaOfCycles, HaessigWanDivisors}]
If $A^1(X)$ is of rank one then all four conjectures are true.
\end{theorem}

Examples of such varieties are nonsingular complete
intersections of dimension 3 or greater and Grassman varieties.
When the rank is greater than one, little is known. However, a
positive result of Wan's states:

\begin{theorem}\label{T: Wan}(Wan \cite{WanZetaOfCycles})
Suppose the effective monoid $A^+(X)$ is finitely generated
over $\bb Z_{\geq 0}$. Let $\rho$ denote the rank of $A^1(X)$
and $h$ the order of its torsion subgroup. Then
\begin{enumerate}
\item $Z_{div}(X,T)$ is $p$-adic meromorphic on the closed
    unit disk $|T|_p \leq 1$.
\item Conjecture II holds. That is, $Z_{div}(X,T)$ has a
    pole of order $\rho$ at $T=1$.
\item The special value at $T=1$ is given by
\[
\lim_{T \rightarrow 1} (1 - T)^\rho Z_{div}(X,T) = \frac{1}{q^\lambda - 1} h R(X)
\]
where there are $q^\lambda-1$ roots of unity in the
function field of $X$, and $R(X)$ is a type of regulator
(see \cite{WanZetaOfCycles} for the definition).
\end{enumerate}
\end{theorem}

When the rank of $A^1(X)$ is strictly greater than one, only
one example currently exists in the literature where all four
conjectures have been proven. This is the quadric surface
$xw=yz$ in $\bb P^3$ which has rank 2
\cite{HaessigWanDivisors}.

All positive results on the zeta function of divisors so far
have been tied with various positive results on the
(generalized) Riemann-Roch problem. For an effective divisor
$D$ on $X$, define $l(D) := dim_{\bb F_q} H^0(X, \cc O(D))$.
The Riemann-Roch problem asks about the behavior of $l(nD)$ as
$n$ varies. For instance, if $D$ is ample then there exists a
polynomial which when evaluated at $n$ agrees with $l(nD)$.
Often, not one but a finite number of polynomials describe the
Riemann-Roch problem. Define a \emph{quasi-polynomial} as a
function $f: \bb Z \rightarrow \bb R$ together with $d$
polynomials $p_1(x), \ldots, p_d(x)$ which satisfy $f(n) =
p_i(n)$ if $n \equiv i$ mod $d$. For projective nonsingular
surfaces defined over $\bb F_q$, $l(nD)$ equals a
quasi-polynomial for all $n$ sufficiently large (see section
\ref{S: RiemannRoch}). This is also true for toric varieties
(see section \ref{S: Toric}). For the zeta function of
divisors, we will need to consider the following
generalization.

\begin{definition}
We say an effective divisor $D$ has \emph{quasi-polynomial
growth} if for each effective divisor $E$ on $X$, $l(E + nD)$
agrees with a quasi-polynomial for all $n$ sufficiently large.
\end{definition}

The main theorem of this paper is the following.

\begin{theorem}\label{T: MainMero}
Let $X$ be a projective variety defined over a finite field.
Assume the effective monoid $A^+(X)$ is finitely generated and
each generator has quasi-polynomial growth. Then
\begin{enumerate}
\item $Z_{div}(X,T)$ is $p$-adic meromorphic.
\item The poles of $Z_{div}(X,T)$ are algebraic over $\bb
    Q_p$ of bounded degree.
\end{enumerate}
\end{theorem}

As mentioned above, when the dimension of $X$ is greater than
one, the zeta function of divisors has zero radius of
convergence over the complex numbers. This suggests that a
weight should be used to compensate, which is what Atsuchi
Moriwaki considers in \cite{MoriwakiZetaCycles}. First, notice
that
\[
Z_{div}(X,T) = \sum_{d = 0}^\infty M_d \> T^d
\]
where $M_d$ is the number of effective divisors on $X$ of
degree $d$. Moriwaki has shown that the series
\[
\sum_{d = 0}^\infty M_d \> T^{d^{dim(X)}}
\]
\emph{does} converge on a small disk around the origin in the
complex numbers. It would be interesting to see what further
properties this series possesses.

\bigskip\noindent{\bf Height zeta function.}
Let $X$ be a (integral) projective variety defined over $\bb
F_q$, and denote its function field by $\bb F_q(X)$. With an
$\bb F_q(X)$-rational point $y := [y_0 : y_1 : \cdots : y_n]
\in \bb P^n(\bb F_q(X))$ we may define the (logarithmic) height
\[
h_X(y) := - deg \inf_i (y_i)
\]
where $\inf_i (y_i)$ denotes the greatest divisor $D$ of $X$
such that $D \leq (y_i)$ for all $i$ with $y_i \not= 0$. Next,
let $Y \hookrightarrow \bb P^n_{\bb F_q(X)}$ be a projective
variety defined over $\bb F_q(X)$. Define the height zeta
function of $Y / \bb F_q(X)$ as
\[
Z_\text{ht}(Y(\bb F_q(X)), T) := \sum_{d=0}^\infty N_d(Y) T^d
\]
where
\[
N_d(Y) := \# \{ y \in Y(\bb F_q(X)) | h_{X}(y) = d \}.
\]

Assuming $A^1(X)$ is of rank one, Wan \cite{WanHeight} proved
$Z_\text{ht}(\bb P^n(\bb F_q(X)), T)$ is always $p$-adic
meromorphic; further, it is rational if and only if $dim(X) =
1$. He conjectured that $p$-adic meromorphic continuation
should hold for any projective variety $X$ defined over $\bb
F_q$ whose effective monoid is finitely generated. From the
methods in this paper and \cite{HaessigHeight}, we are able to
support his conjecture by the following:

\begin{theorem}
Let $X$ be a projective variety defined over a finite field.
Assume the effective monoid $A^+(X)$ is finitely generated and
each generator has quasi-polynomial growth. Then
$Z_\text{ht}(\bb P^n(\bb F_q(X)), T)$ is $p$-adic meromorphic.
\end{theorem}

\bigskip\noindent{\bf Acknowledgements.} I would like to thank
Mircea Musta\c{t}\v{a} for helpful comments about toric
varieties, and Vasudevan Srinivas and Dale Cutkosky for
insightful discussions on surfaces over finite fields. Lastly,
I would like to thank the referee for their very helpful
suggestions.

\section{Riemann-Roch problem for surfaces}\label{S: RiemannRoch}

The zeta function of divisors was initially defined and studied
by Wan \cite{WanZetaOfCycles}. In that paper he establishes a
method of study based on the Riemann-Roch approach of F. K.
Schmidt \cite{Schmidt} (See \cite{Roquette} for an historical account)
who used this method to prove rationality of the Weil zeta function of curves defined over a finite field. An initial obstruction to this line of
investigation is an understanding of the so-called Riemann-Roch
problem.

Given a divisor $D$, the Riemann-Roch problem is concerned with
the dimension variation of the space of global sections $H^0(X,
\cc O(nD))$ as $n$ varies. This variation is often very well
behaved. Indeed, when $D$ is ample, the dimension $l(nD) := dim
H^0(X, \cc O(D))$ equals a (Hilbert-Samuel) polynomial of
degree $dim(X)$ evaluated at $n$ for all $n$ sufficiently
large. When this type of behavior is expected in general seems
to be an open question. However, it certainly does not always
occur as an example of Cutkosky and Srinivas demonstrates
\cite[Section 7]{CutkoskyRRproblem}. They describe an effective
divisor $D$ on a nonsingular projective 3-fold defined over a
finite field for which $l(nD)$ is a polynomial of degree 3 in
$\left\lfloor n(2- \sqrt{3}/3) \right\rfloor$.

The quasi-polynomial behavior for surfaces discussed below is a
consequence of a deeper algebraic result. For each divisor $D$
define the graded $\bb F_q$-algebra
\[
R[X, D] := \bigoplus_{n \geq 0} H^0(X, \cc O(nD)).
\]
Next, for each divisor $E$ on $X$ define the $R[X, D]$-module
\[
R[X, D; E] := \bigoplus_{n \geq 0} H^0(X, \cc O( E + nD )).
\]
If $R[X,D]$ is finitely generated then for all $n$ sufficiently
large $l(nD)$ equals a quasi-polynomial in $n$. Furthermore,
when $R[X,D]$ is finitely generated and $R[X,D;E]$ is finite
over $R[X,D]$ then $l(E + nD)$ equals a quasi-polynomial in
$n$.

\begin{question}\label{Q: finite}
When is $R[X,D]$ finitely generated, and when is $R[X, D; E]$ a
finitely generated graded $R[X,D]$-module?
\end{question}

For $D$ a semi-ample divisor then a positive answer to the
latter part of Question \ref{Q: finite} is known. Recall, $D$
is a semi-ample divisor if there exists a positive integer $m$
such that the complete linear system $|mD|$ is base point free.
In this case, for any projective variety $X$ defined over an
arbitrary field, Zariski \cite[Theorem 5.1]{ZariskiRR} has
shown that for every effective divisor $E$ and semi-ample
divisor $D$, $R[X, D; E]$ is a finitely generated graded $R[X,
D]$-module.

\begin{theorem}
Let $X$ be a projective nonsingular surface defined over the
finite field $\bb F_q$. Let $D$ be an effective divisor on $X$.
Then,
\begin{enumerate}
\item (Zariski \cite{ZariskiRR}, Cutkosky and Srinivas
    \cite{CutkoskyRRproblem}) $R[X,D]$ is a finitely
    generated, graded $\bb F_q$-algebra.
\item $R[X, D; E]$ is a finite $R[X, D]$-module for every
    effective divisor $E$ on $X$.
\end{enumerate}
\end{theorem}

\begin{proof}
That $R[X,D]$ is finitely generated was proven by Zariski
\cite{ZariskiRR} and Cutkosky and Srinivas
\cite{CutkoskyRRproblem}. In order to prove the second part of
the theorem, we will need to recall their methods. Let $D$ be
an effective divisor of $X$, and denote by $P_1, \ldots, P_m$
its prime components. We may define a (symmetric) quadratic
form $\phi_D$ by $\sum (P_i \cdot P_j) x_i x_j$. The
eigenvalues associated to this form are real, and all but at
most one can be positive \cite[p.588]{ZariskiRR}. We say $D$ is
of type $(s,t)$ if the associated form has $s$ positive
eigenvalues and $t$ negative eigenvalues; that is, $D$ may be
of type $(1,t)$ or $(0,t)$.

\bigskip\noindent\emph{Case $\phi_D$ is of type $(0,t)$.} In
\cite[Theorem 11.5]{ZariskiRR}, Zariski proves $R[X, D]$ is a
finitely generated $\bb F_q$-algebra. Using his notation, for
some $n', e \in \bb Z_{>0}$, he demonstrates that the complete
linear system $|n' e D|$ has no fixed components. By
\cite[Theorem 6.1]{ZariskiRR}, this means $|m n' eD|$ is base
point free for all $m$ sufficiently large. It follows from
\cite[Theorem 5.1]{ZariskiRR} that $R[X, D; E]$ is a finite
$R[X,D]$-module.

\bigskip\noindent\emph{Case $\phi_D$ is of type $(1,t)$.}
By Zariski Decomposition \cite[Theorem 7.7]{ZariskiRR} (see
\cite{BauerQuickZariskiDecomp} for a quick proof) there exists
a unique $\bb Q$-divisor $\cc E$, called the arithmetically
negative part of $D$, for which by \cite[Theorem
10.6]{ZariskiRR} $R[X,D]$ is a finitely generated $\bb
F_q$-algebra if and only if $|m_0(D - \cc E)|$ has no fixed
components for some $m_0$ (hence $D- \cc E$ is semi-ample).
Note, $m_0 \cc E$ is a $\bb Z$-divisor. In the proof of this,
Zariski says $n m_0 \cc E$ is a fixed component of $|n m_0 D|$
for all $n$, which follows from \cite[Theorem 8.1]{ZariskiRR}.
Thus, we have an equality between linear systems:
\[
\{ f \in \bb F_q(X) | (f) + nm_0 D \geq 0\} = \{ f \in \bb F_q(X) | (f) + nm_0 D - n m_0 \cc E \geq 0\}.
\]
Consequently,
\[
R[X, m_0 D] = R[X, m_0(D - \cc E)].
\]
Zariski was unable to prove $D-\cc E$ is semi-ample. However,
Cutkosky and Srinivas in the proof of \cite[Theorem
3]{CutkoskyRRproblem} indeed show $D- \cc E$ is semi-ample.
Thus, since
\begin{equation}\label{E: ZariskiRings}
R[X,D] = \bigoplus_{r=0}^{m_0-1} R[X, m_0 D ; rD] = \bigoplus_{r=0}^{m_0-1} R[X, m_0 (D-\cc E) ; rD]
\end{equation}
and $R[X, m_0(D - \cc E); rD]$ is a finite $R[X, m_0(D-\cc
E)]$-module (hence a finite $R[X, m_0 D]$-module), we see that
$R[X, D]$ is a finite $R[X, m_0 D]$-module. Since $R[X,D]$ is
an algebra, and $R[X, m_0 D]$ is a finitely generated algebra,
it follows that $R[X, D]$ is a finitely generated $\bb
F_q$-algebra.

A similar argument follows for $R[X,D; E]$: decompose $R[X,
D;E]$ as in (\ref{E: ZariskiRings}). Then $R[X, m_0(D - \cc E);
rD + E]$ is a finite $R[X, m_0(D - \cc E)]$-module, hence a
finite $R[X, m_0 D]$-module, and thus, a finite $R[X,
D]$-module.
\end{proof}

\begin{corollary}
Let $X$ be a projective nonsingular surface defined over the
finite field $\bb F_q$. All effective divisors $D$ on $X$ have
quasi-polynomial growth.
\end{corollary}

While divisors on a nonsingular surface defined over a finite
field may have quasi-polynomial growth, it may happen that the
effective monoid $A^+(X)$ is not finitely generated. For
instance, it is known that on a nonsingular projective rational
surface whose anticanonical divisor $-K_X$ is effective,
$A^+(X)$ is finitely generated if and only if $X$ has only a
finite number of $(-1)$-curves and only a finite number of
$(-2)$-curves; examples where such curves are finite in number
may be found in \cite{Lahyane}. However, this cannot happen on
toric varieties.

In the next section we will see that $A^+(X)$ is always
finitely generated for projective toric varieties. Furthermore,
every effective divisor $D$ has quasi-polynomial growth.

\section{Riemann-Roch problem for Toric Varieties}\label{S: Toric}

For toric varieties, the Riemann-Roch problem has an elegant
solution related to counting integer points in polytopes. For
any T-invariant Weil divisor $D$ on a complete toric variety,
$l(nD)$ equals the number of lattice points in an associated
polytope $P_{nD}$. Since $P_{nD} = n P_D$, we see that the
Riemann-Roch problem is equivalent to studying the number of
lattice points in a dilation of the polytope $P_D$. This was
extensively studied, and answered, by Ehrhart. In particular,
$l(nD)$ is again a quasi-polynomial in $n$. Let us be more
precise.

We will freely use results and terminology from Fulton
\cite{FultonToric}. While Fulton's results are mostly stated
over the complex numbers, only an algebraically closed field is
needed for the results we will use.

Let $k$ be an algebraically closed field. Let $X$ be a
$d$-dimensional projective toric variety over $k$, defined by a
complete fan $\Delta$ in a lattice $N \cong \bb Z^d$. Let $M :=
Hom(N, \bb Z)$ be the dual lattice of $N$, and denote by
$\langle \cdot, \cdot \rangle$ the bilinear pairing on $M$ and
$N$. Let $\Delta(1)$ be the set of one-dimensional rays (cones)
in $\Delta$. For each $\rho \in \Delta(1)$, denote by $v_{\rho}
\in N$ the unique generator of $\rho \cap N$.

There is an action of the $d$-torus $T := N \otimes_{\bb Z} k^*
= Hom_{\bb Z}(N, k^*)$ on $X$. For every $\rho \in \Delta(1)$
there corresponds an irreducible $T$-invariant Weil divisor
$D_{\rho}$ which is the orbit closure of $T$ acting on the
distinguished point associated to $\rho$. Further, every
$T$-invariant Weil divisor is a linear combination of these.

It follows from the following exact sequence that the group
$A^1(X)$ is completely determined by the torus action:
\begin{equation}\label{E: exactSeq}
0 \rightarrow M \rightarrow \bigoplus_{\rho \in \Delta(1)} \bb Z
D_{\rho} \rightarrow A^1(X) \rightarrow 0.
\end{equation}
where the second map is $m \mapsto \sum_{\rho \in \Delta(1)}
\langle m, v_\rho \rangle D_{\rho}$. This exact sequence also
helps determine the effective monoid $A^+(X)$ as Theorem
\ref{T: ToricEffec} below demonstrates.

Associated to each $T$-invariant Weil divisor $D = \sum_{\rho
\in \Delta(1)} b_\rho D_\rho$ is a polyhedron $P_D$ defined by:
\begin{equation} \label{E: Poly}
P_D := \{ m \in M \otimes_{\bb Z} \bb R : \langle m,
v_{\rho} \rangle \geq -b_{\rho} \text{ for each } \rho \}.
\end{equation}
This polyhedron is compact since the fan $\Delta$ is complete.
There is a bijection between a basis of the Riemann-Roch space
$H^0(X, \cc O(D))$ and the number of lattice points in the
polyhedron $P_D$:
\begin{equation}\label{E: RRPolytope}
l(D) = \#(P_D \cap M).
\end{equation}

\begin{theorem}\label{T: ToricEffec}
For $X$ a projective toric variety, the effective monoid
$A^+(X)$ is finitely generated. In particular, it is generated
by the image of $\bigoplus_{\rho \in \Delta(1)} \bb Z_{\geq 0}
D_{\rho}$ in $A^1(X)$.
\end{theorem}

\begin{proof}
Let $E$ be an effective divisor on $X$. By (\ref{E: exactSeq}),
$E$ is linearly equivalent to a $T$-invariant divisor $D$, thus
$l(D) = \# (P_D \cap M) \geq 1$. Next, write $D = \sum_{\rho
\in \Delta(1)} b_\rho D_\rho$. Let $A$ be the matrix with rows
the vectors $v_{\rho}$ and let $b$ be the column vector with
entries $b_{\rho}$. We say $v \geq 0$ if its entries are each
greater than or equal to zero. With this new notation, and
viewing the pairing between $N$ and $M$ as an inner product, we
may rewrite (\ref{E: Poly}) as
\[
P_D = \{ m \in M \otimes_{\bb Z} \bb R : Am \geq -b \}.
\]
where we represent $D$ as the vector $b$. Since $P_D \cap M$ is
non-empty, there exists vectors $m \in M$ and $q \geq 0$ such
that $Am + b = q$. In terms of divisors, this means $D$ is
linearly equivalent to the divisor $D_q$ defined by $q$. The
theorem follows since $D_q$ is an element of the effective
divisors $\bigoplus_{\rho \in \Delta(1)} \bb Z_{\geq 0}
D_{\rho}$.
\end{proof}

\begin{theorem}\label{T: ToricRR}
Let $E$ and $D$ be effective $T$-invariant divisors on a
projective toric variety of dimension $d$. Then $l(E + nD)$
equals a quasi-polynomial, whose degree is bounded below by 1
and above by $d$, for all nonnegative integers $n$.
\end{theorem}

\begin{proof}
From (\ref{E: RRPolytope}) we have $l(E + nD) = \# (P_{nD + E}
\cap M)$ for every $n$. Write $D = \sum_{\rho \in \Delta(1)}
b_\rho D_\rho$ and $E = \sum_{\rho \in \Delta(1)} a_\rho
D_\rho$ with $b_\rho, a_\rho \in \bb Z_{\geq 0}$. As in the
proof of Theorem \ref{T: ToricEffec}, let $A$ be the matrix
with rows the vectors $v_{\rho}$, let $b$ be the column vector
with entries $b_{\rho}$, and let $a$ be the column vector with
entries $a_{\rho}$. Then
\[
P_{nD+E} = \{ m \in M \otimes_{\bb Z} \bb R : Am \geq -nb - a \}.
\]
Using terminology from \cite{ClaussEhrhart}, since $P_{nD} = n
P_D$, the polytope $P_{E + nD}$ constitutes a \emph{bordered
system}. By \cite[Theorem 2]{ClaussEhrhart}, the number of
lattice points in $P_{nD+E}$ equals a quasi-polynomial of
degree bounded above by $d$. Furthermore, this holds for all
$n$. Since the dimension of $P_{nD}$ is at least one, the
degree of the quasi-polynomial is at least one.
\end{proof}

While Theorem \ref{T: ToricRR} allows us to prove $p$-adic
meromorphic continuation for the zeta function of divisors, the
proof is rather unsatisfying since the author believes the
stronger result is true, that $R[X,D;E]$ is finite over
$R[X,D]$. Note, Elizondo \cite{ElizondoMutiples} has already
shown that $R[X,D]$ is finitely generated for any effective
divisor $D$.

\section{$p$-adic Meromorphic continuation}

In this section, we will prove a $p$-adic meromorphic
continuation theorem for the zeta function of divisors which
applies whenever the generators of $A^+(X)$ have
quasi-polynomial growth. The proof will also show that
adjoining all the poles of the zeta function to $\bb Q_p$
creates a finite extension field.

When the effective monoid $A^+(X)$ is finitely generated, Wan
\cite{WanZetaOfCycles} reduces the zeta function of divisors to a finite sum of both rational functions, and series of the form
\[
Z_E(X,T) := \sum_{n_1, \ldots, n_r \geq 0} q^{l(E + n_1 D_1 + \cdots n_r D_r)}
T^{e + n_1 d_1 + \cdots + n_r d_r}
\]
where $E, D_1, D_2, \ldots, D_r$ are effective divisors, $e :=
deg(E)$, and $d_i := deg(D_i)$. Thus, we need only study the
$p$-adic analytic behavior of these series.

\begin{lemma}\label{L: ExtendedRR}
Let $X$ be a projective variety. Let $D_1, \ldots, D_r$ be
effective divisors on $X$, each having quasi-polynomial growth.
Then for every effective divisor $E$ on $X$ there exists a
positive integer $N$ such that $l(E + n_1 D_1 + \cdots + n_r
D_r)$ equals a quasi-polynomial in the variables $n_1, \ldots,
n_r$ for all $n_i \geq N$.
\end{lemma}

\begin{proof}
This was first proven in \cite[Theorem 2]{ClaussEhrhart} in the
context of counting lattice points in polytopes. For
completeness, we will present it here. The proof follows by
induction on $r$. Assume the result for $r-1$ divisors. Fixing
positive integers $n_1, \ldots, n_{r-1}$ and letting $n_r$
vary, since $D_r$ has quasi-polynomial growth there exists
$N_r$ such that
\[
l(E + n_1 D_1 + \cdots + n_r D_r) = \alpha_d n_r^d + \alpha_{d-1}
n_r^{d-1} + \cdots + \alpha_0 \qquad \text{for all } n_r \geq
N_r,
\]
where the right-hand side is a quasi-polynomial. Note, the
coefficients $\alpha_i$ depend on the fixed $n_1, \ldots,
n_{r-1}$. Plugging in values for $n_r$, we obtain a linear
system
\[
\begin{cases}
\alpha_d (N_r+0)^d + \alpha_{d-1}(N_r+0)^{d-1} + \cdots + \alpha_0
&=
l(E + n_1 D_1 + \cdots + (N_r+0) D_r) \\
\alpha_d (N_r+1)^d + \alpha_{d-1}(N_r+1)^{d-1} + \cdots + \alpha_0
& = l(E + n_1 D_1 + \cdots + (N_r+1) D_r) \\
\qquad \qquad \qquad \vdots & \qquad \qquad \vdots \\
\end{cases}
\]
By Cramer's rule, each coefficient $\alpha_i$ is a linear
combination of elements from the right-hand side. Since, by the
induction hypothesis, each member of the right-hand side is a
quasi-polynomial for all large $n_1, \ldots, n_{r-1}$, we see
that $\alpha_i$ will also be of this form.
\end{proof}

From this lemma and the discussion at the beginning of this
section, the zeta function of divisors may be reduced to a
finite sum of both rational functions, and series of the form
\[
\sum_{x_1, \ldots, x_n \geq 0} q^{l(x_1, \ldots, x_n)} T^{d_1 x_1
+ \cdots + d_n x_n}
\]
where $l(x_1, \ldots, x_n)$ is an increasing polynomial with
rational coefficients in the variables $x_1, \ldots, x_n$.

\begin{theorem}\label{T: MeroResult}
Let $f \in \bb Q[x_1, \ldots, x_n]$ be an increasing
polynomial; that is,
\[
f(x_1, \ldots, x_i + 1, \ldots, x_n) \geq f(x_1, \ldots, x_i,
\ldots, x_n) \quad \text{for all } x_1, \ldots, x_n \in \bb Z_{\geq 0}.
\]
Let $d_1, \ldots, d_n$ be positive integers. Then
\[
F(T) := \sum_{x_1, \ldots, x_n \geq 0} q^{f(x_1, \dots, x_n)} T^{d_1
x_1 + \cdots + d_n x_n}
\]
is $p$-adic meromorphic on $\bb C_p$. Furthermore, adjoining
the poles to $\bb Q_p$ creates a finite extension field of $\bb
Q_p$.
\end{theorem}

The proof will consist of the rest of this section and will
follow by induction on the number of variables $n$. In the case
when there is only one variable, $n = 1$, the result follows by
either an application of the geometric series (when $f$ is a
linear polynomial) in which case $F(T)$ is meromorphic, else it
is entire (when $f$ has degree at least two).

Suppose the result is true for $1, 2, \ldots, n-1$ variables. Let $S_n$
denote the symmetric group on $\{ 1, \ldots, n\}$. A
permutation $\sigma \in S_n$ acts on points by permuting the
coordinates:
\[
\sigma(x_1, \ldots, x_n) := (x_{\sigma(1)}, \ldots, x_{\sigma(n)}).
\]
Define
\[
E^{(n)} := \{ (k_1 + k_2 + \cdots + k_n, k_2 + k_3 + \cdots + k_n,
\ldots, k_n) | k_1, \ldots, k_n \in \bb Z_{\geq 0} \}.
\]
For each $\sigma \in S_n$, define the set $E_\sigma^{(n)} :=
\sigma(E^{(n)})$. Notice that $\bb Z_{\geq 0}^n =
\bigcup_{\sigma \in S_n} E_\sigma^{(n)}$, and that points in $\bb
Z_{\geq 0}^n$ may be in several $E_\sigma^{(n)}$. We would like
to use the inclusion-exclusion principle on this collection. To
accomplish this, let $\cc E_n$ denote the collection of all
possible intersections between elements of $\{ E_\sigma^{(n)}
\}_{\sigma \in S_n}$ (that is, $\cc E_n$ is the closure of the
set $\{ E_\sigma^{(n)} \}_{\sigma \in S_n}$ by the operation of
intersection). Observe that elements in $\cc E_n$ are of the
form:
\[
\{ \sigma(k_1 + k_2 + \cdots + k_r, \ldots, k_r) | k_1, \ldots, k_r \in \bb Z_{\geq 0} \},
\]
for some $1 \leq r \leq n$ and $\sigma \in S_n$ with possible
repetition of the coordinates in the middle. For example, with
$n=3$ and $r =2$, $\cc E_3$ contains a set of the form
\[
\{ (k_1 + k_2, k_1+k_2, k_2)  | k_1, k_2 \in \bb Z_{\geq
0} \}.
\]
Consequently, using the inclusion-exclusion principle, we may write
\begin{equation}\label{E: inductionSum}
F(T) = \sum_{H \in \cc E_n} c_H \sum_{(x_1, \ldots, x_n) \in H}
q^{f(x_1, \ldots, x_n)} T^{d_1 x_1 + \cdots d_n x_n}.
\end{equation}
where $c_H$ is some constant coming from the
inclusion-exclusion. Since $\cc E_n$ is a finite set, we are
reduced to proving $\sum_{(x_1, \ldots, x_n) \in H} q^{f(x_1,
\ldots, x_n)} T^{d_1 x_1 + \cdots d_n x_n}$ is $p$-adic
meromorphic for each $H \in \cc E_n$.

For a typical sum in (\ref{E: inductionSum}), we have
\begin{equation}\label{E: inductionSum2}
\sum_{k_1, \ldots, k_r \geq 0} q^{f(k_1 + \cdots + k_r,
\ldots, k_r)} T^{d_1(k_1 + \cdots + k_r) + \cdots + d_n(k_r)}
\end{equation}
Notice that, unless $H$ equals $E_\sigma^{(n)}$ for some
$\sigma$, we may write
\[
f(k_1 + \cdots + k_r, \ldots, k_r) = g(k_1, \ldots, k_r)
\]
where $1 \leq r \leq n-1$ and $g$ is an increasing polynomial
in the variables $k_1, \ldots, k_r$. It follows by the
induction hypothesis that (\ref{E: inductionSum2}) is $p$-adic
meromorphic. Thus, we need only consider the case when $H =
E_\sigma^{(n)}$.

Without loss of generality, suppose $H = E^{(n)}$. Thus, with
\begin{equation}\label{E: subs}
\begin{split}
x_1 &= k_1 + k_2 + k_3 + \cdots + k_n \\
x_2 &= \phantom{k_1 +} k_2 + k_3 + \cdots + k_n \\
&\vdots \qquad\qquad\qquad \ddots \\
x_n &= \phantom{k_1 + k_2 + k_3 + \cdots +} k_n,
\end{split}
\end{equation}
let us prove
\begin{equation}\label{E: F1}
F_1(T) := \sum_{k_1, \ldots, k_n \geq 0} q^{f(k_1 + \cdots + k_n,
\ldots, k_n)} T^{d_1(k_1 + \cdots + k_n) + \cdots + d_n(k_n)}
\end{equation}
is $p$-adic meromorphic. Our argument breaks up into two parts:
when $deg_{x_1}(f) \geq 2$ and when $deg_{x_1}(f) = 1$.

Let us first consider the case when $deg_{x_1}(f) \geq 2$. We
need the following technical lemma.

\begin{lemma}\label{L: technicalLemma}
For each $i = 1, \ldots, n$, write
\[
f(k_1 + \cdots + k_n, \ldots, k_n) = A_{r_i}^{(i)} k_i^{r_i} + A_{r_i-1}^{(i)} k_i^{r_i-1} + \cdots + A_0^{(i)}.
\]
Then $A_{r_i}^{(i)}$ is an increasing polynomial which depends
only on $k_{i+1}, \ldots, k_n$. Furthermore, $A_{r_i}^{(i)}$ is
always nonnegative, and $r_i \geq \max_{1 \leq j \leq i}
deg_{x_j}(f)$ for each $i$.
\end{lemma}

\begin{proof}
Notice that $f$ is an increasing polynomial in $k_1, \ldots,
k_n$. It is not hard to show that the coefficient
$A_{r_i}^{(i)}$ of the largest power of $k_i$ is a polynomial
in $k_{i+1}, \ldots, k_n$. Let us show it is an increasing
polynomial. Write $A_{r_i}^{(i)} = g_{r_i}(k_{i+1}, \ldots,
k_n)$. Writing $f(k_1, \ldots, k_n)$ for $f(k_1 + \cdots + k_n,
\ldots, k_n)$, then for any $j > i$ we have
\begin{align*}
f(k_1, \ldots, k_j + 1, &\ldots, k_n) - f(k_1, \ldots, k_n) \\
&= \bigg( g_{r_i}(k_{i+1}, \ldots, k_j + 1, \ldots, k_n) - g_{r_i}(k_{i+1},
\ldots, k_n) \bigg) k_i^{r_i} + (\text{lower terms}).
\end{align*}
Since $f$ is an increasing function, this must be nonnegative
for all $k_1, \ldots, k_n \geq 0$, and so the coefficient of
$k_i^{r_i}$ must be nonnegative. This shows $g_{r_i}$ is an
increasing polynomial.

Next, let us show $r_i \geq \max_{1 \leq j \leq i}
deg_{x_j}(f)$ and $g_{r_i}$ is not identically zero. Write
\[
f(x_1, \ldots, x_n) = \sum_{j_1, \ldots, j_n} C_{j_1, \ldots, j_n} x_1^{j_1} \cdots x_n^{j_n}
\]
where the $C_{j_1, \ldots, j_n}$ are nonzero rational numbers.
Let $R_i$ be the maximum of all $j_1 + \cdots + j_i$. We will
show $R_i = r_i$. Write
\begin{equation}\label{E: technical2}
f(x_1, \ldots, x_n) = \sum_{j_1 + \cdots + j_i = R_i} H_{j_1, \ldots, j_i}(x_{i+1}, \ldots, x_n)
x_1^{j_1} \cdots x_i^{j_i} + (\text{lower terms})
\end{equation}
where $H_{j_1, \ldots, j_i}$ is a polynomial in $x_{i+1},
\ldots, x_n$. Note, ``lower terms'' means $x_1^{l_1} \cdots
x_n^{l_n}$ with $l_1 + \cdots + l_n < R_i$. Using the
substitution (\ref{E: subs}), we have
\begin{equation}\label{E: technical}
f(k_1, \ldots, k_n) = \sum_{j_1 + \cdots + j_i = R_i} H_{j_1, \ldots, j_i}(k_{i+1} + \cdots + k_n,
\ldots, k_n) (k_1 + \cdots + k_n)^{j_1} \cdots (k_i
+ \cdots + k_n)^{j_i} + (\text{lower terms}).
\end{equation}
Notice that the largest power of $k_i$ to appear is
$k_i^{R_i}$, and furthermore, by the definition of $R_i$, must
satisfy $R_i \geq \max_{1 \leq j \leq i} deg_{x_j}(f)$. Let us
show the coefficient of $k_i^{R_i}$ is nonzero (and hence $r_i
= R_i$).

The coefficient of $k_i^{R_i}$ in (\ref{E: technical}) is
\[
\sum_{j_1 + \cdots + j_i = R_i} H_{j_1, \ldots, j_i}(k_{i+1} + \cdots + k_n,
\ldots, k_n),
\]
and this contains the unique polynomial
\[
\sum_{j_1 + \cdots + j_i = R_i} H_{j_1, \ldots, j_i}(k_{i+1}, \ldots, k_n)
\]
which is a nonzero polynomial by (\ref{E: technical2}).
\end{proof}

Let us show $F_1$ is $p$-adic meromorphic when $deg_{x_1}(f)
\geq 2$. By Lemma \ref{L: technicalLemma}, the coefficient of
the largest power of $k_n$ is a positive rational number, and
since the degree of $x_1$ is at least 2, the degree of $k_n$ is
also at least two. Consider the coefficient of the largest
power of $k_{n-1}$ in $f$. By Lemma \ref{L: technicalLemma},
this is a nonzero increasing polynomial in $k_n$. Find $k_n'$
such that this polynomial is positive for all $k_n \geq k_n'$.
Then
\[
F_1(T) = \left( \sum_{\substack{k_1, \ldots, k_{n-1} \geq 0 \\
0 \leq k_n < k_n'}} + \sum_{\substack{k_1, \ldots, k_{n-1}
\geq 0 \\ k_n \geq k_n'}} \right) q^{f(k_1 + \cdots + k_n, \ldots, k_n)} T^{d_1 k_1 + \cdots + d_n k_n}.
\]
Notice that the first sum is $p$-adic meromorphic by the
induction hypothesis (it has fewer than $n$ variables). Thus,
let us proceed with the second summation.

By Lemma \ref{L: technicalLemma}, the coefficient of the
largest power of $k_{n-2}$ is a nonzero increasing polynomial
in $k_{n-1}$ and $k_n$. Find constants $k_{n-1}'$ and $k_n''$
such that this polynomial is always positive whenever $k_{n-1}
\geq k_{n-1}'$ and $k_n \geq k_n''$. Then we may write
\[
(\text{second sum}) = \left( \sum_{\substack{k_1, \ldots, k_n \geq 0 \\
0 \leq k_{n-1} < k_{n-1}' \text{ or } k_n' \leq k_n < k_n''}} + \sum_{\substack{k_1, \ldots, k_{n-2}
\geq 0 \\ k_{n-1} \geq k_{n-1}' \text{ and } k_n \geq k_n''}} \right) q^{f(k_1 + \cdots + k_n, \ldots, k_n)} T^{d_1 k_1 + \cdots + d_n k_n}.
\]
As before, the first summation is $p$-adic meromorphic by the
induction hypothesis.

Continuing this argument, we are finally left with a series of
the form
\[
\sum_{k_1, \ldots, k_n \gg 0} q^{f(k_1 + \cdots + k_n,
\ldots, k_n)} T^{d_1(k_1 + \cdots + k_n) + \cdots + d_n(k_n)}
\]
where $deg_{k_i}(f_i) \geq 2$ and the coefficient of the
largest power of $k_i$ is always positive. It follows from the
proof of Lemma \ref{L: MeroLemma} below that this is a $p$-adic
entire function. This finishes the proof of the result when
$deg_{x_1}(f) \geq 2$.

\bigskip

Consider now the case when $deg_{x_1}(f) = 1$. In this case, if
we write
\[
f(x_1, \ldots, x_n) = g_1(x_2, \ldots, x_n) x_1 + g_0(x_2, \ldots,
x_n),
\]
then by an application of geometric series,
\[
F_1(T) = \sum_{k_2, \ldots, k_n \geq 0} \frac{q^{g_1(k_2 + \cdots +
k_n, \ldots, k_n) ( k_2 + \cdots + k_n ) + g_0( k_2 + \cdots + k_n,
\ldots, k_n) } T^{d_2' k_2 + \cdots + d_n' k_n}}{1 - q^{g_1(k_2 +
\cdots + k_n, \ldots, k_n)} T^{d_1}}
\]
where $d_2' := d_1 + d_2$, $d_3' := d_1 + d_2 + d_3$, etc. This
will be a meromorphic function by the following lemma. (Note:
$g_0$ is an increasing polynomial by taking $x_1 = 0$.)

\begin{lemma}
Let $g_1$ and $g_0$ be increasing polynomials. Then
\[
F_1(T) := \sum_{k_2, \ldots, k_n \geq 0} \frac{q^{g_1(k_2 + \cdots + k_n,
\ldots, k_n) ( k_2 + \cdots + k_n ) + g_0( k_2 + \cdots + k_n,
\ldots, k_n) } T^{d_2 k_2 + \cdots + d_n k_n}}{1 - q^{g_1(k_2
+ \cdots + k_n, \ldots, k_n)} T^{d_1}}
\]
is $p$-adic meromorphic.
\end{lemma}

\begin{proof}
We will proceed by induction on the number of variables.
Suppose $n=2$, then $F_1$ takes the form
\[
\sum_{k_2 \geq 0} \frac{q^{g_1(k_2) k_2 + g_0(k_2)} T^{d_2 k_2}}{1-q^{g_1(k_2)} T^{d_1}}.
\]
If either $deg_{k_2}(g_1) \geq 1$ or $deg_{k_2}(g_0) \geq 2$
then this summation is $p$-adic meromorphic by Lemma \ref{L:
MeroLemma} below (note: we are using the fact that $g_1$ and
$g_0$ are increasing polynomials, and so the coefficients of
the largest powers of $k_2$ are necessarily positive). If $g_1$
is a constant (perhaps zero) and $deg_{k_2}(g_0) = 1$ then the
series is meromorphic by an application of geometric series.

Let us assume the result is true whenever there are fewer than $n$ variables. For convenience, let us define
\begin{align*}
u(k_2, \ldots, k_n) :&= \text{ the exponent of $q$ in the numerator of $F_1$} \\
&= g_1(k_2 + \cdots + k_n,
\ldots, k_n) ( k_2 + \cdots + k_n ) + g_0( k_2 + \cdots + k_n,
\ldots, k_n).
\end{align*}
In the following, we will move progressively backwards through
the variables: $k_n, k_{n-1}, \ldots, k_2$.

By Lemma \ref{L: technicalLemma} we know the degree of $k_n$ in
$u$ is at least one and the coefficient of this largest degree
term is a positive rational number. If the degree is equal to
one, then $g_1$ must be identically zero and $g_0$ is linear in
$k_2, \ldots, k_n$. It follows that $F_1$ is meromorphic.

Suppose the degree of $k_n$ in $u$ is at least two. By Lemma
\ref{L: technicalLemma}, the only way for $k_{n-1}$ to not
appear in $u$ is if $g_1$ is identically zero and $g_0$ depends
only on $k_n$. (Again, we are using the fact that $g_1$ and
$g_0$ are increasing polynomials, and so the coefficient of the
largest power of $k_{n-1}$ in both $g_1$ and $g_0$ must be
nonnegative. In particular, these terms cannot eliminate one
another.) It follows from Lemma \ref{L: MeroLemma} below that
this is a meromorphic function. Thus, let us suppose the degree
of $k_{n-1}$ in $u$ is at least one. By Lemma \ref{L:
technicalLemma} the coefficient of the largest degree term of
$k_{n-1}$ is a nonzero increasing polynomial in $k_n$. Find
$k_n'$ large enough such that for all $k_n \geq k_n'$ this
polynomial takes positive values. From this, we may write
\[
F_1(T) = \left( \sum_{ \substack{k_2, \ldots, k_{n-1} \geq 0 \\
0 \leq k_n < k_n'}} + \sum_{ \substack{k_2, \ldots, k_{n-1} \geq 0 \\ k_n \geq k_n'}} \right)
\frac{q^{u(k_2, \ldots, k_n)} T^{d_2 k_2 + \cdots + d_n k_n}}{1 - q^{g_1(k_2
+ \cdots + k_n, \ldots, k_n)} T^{d_1}}
\]
Notice that the first summation takes on a similar form to
$F_1(T)$ but with variables $k_2, \ldots, k_{n-1}$. Hence, it
is meromorphic by the induction hypothesis. Let us consider the
second summation
\[
F_2(T) := \sum_{ \substack{k_2, \ldots, k_{n-1} \geq 0 \\ k_n \geq k_n'}}
\frac{q^{u(k_2, \ldots, k_n)} T^{d_2 k_2 + \cdots + d_n k_n}}{1 - q^{g_1(k_2
+ \cdots + k_n, \ldots, k_n)} T^{d_1}}.
\]
where the degree of $k_{n-1}$ in $u$ is at least one and the
degree of $k_n$ in $u$ is at least two, and the coefficients of
their largest power is always positive. Suppose the degree of
$k_{n-1}$ is one. Then $g_1$ must be a nonnegative constant
$c_1$ and $g_0$ cannot contain any $k_2, \ldots, k_{n-2}$. Thus
$F_2(T)$ takes the form
\[
\sum_{ \substack{k_2, \ldots, k_{n-1} \geq 0 \\ k_n \geq k_n'}}
\frac{{ q^{c_1(k_2 + \cdots + k_n) + c_2 k_{n-1}
+ h(k_n)}} T^{d_2 k_2 + \cdots + d_n k_n}}{1 - q^{c_1} T^{d_1}}
\]
where $h$ is a quadratic polynomial in $k_n$, and $c_2$ is a
constant. From geometric series and Lemma \ref{L: MeroLemma}
below, this is a meromorphic function. Thus, we may suppose
that the degrees of both $k_{n-1}$ and $k_n$ in $u$ in $F_2(T)$
are at least two, and the coefficients of their largest powers
are always positive.

A similar argument works for $k_{n-2}$ as follows. By Lemma
\ref{L: technicalLemma}, the only way for $k_{n-2}$ to not
appear in $u$ is if $g_1$ is identically zero and $g_0$ is a
polynomial in $k_{n-1}$ and $k_n$. By our assumptions on $u$,
it follows from Lemma \ref{L: MeroLemma} below that this is a
meromorphic function. Thus, let us suppose the degree of
$k_{n-2}$ in $u$ is at least one.

By Lemma \ref{L: technicalLemma} the coefficient of the largest
degree term of $k_{n-2}$ is a nonzero increasing polynomial in
$k_{n-1}$ and $k_n$. Find $k_{n-1}'$ and $k_n''$ large enough
such that whenever both $k_{n-1} \geq k_{n-1}'$ and $k_n \geq
k_n''$ this polynomial takes only positive values. From this,
assuming $k_n' \geq k_n''$, we may write
\[
F_2(T) = \left( \sum_{ \substack{k_2, \ldots, k_n \geq 0 \\
0 \leq k_{n-1} < k_{n-1}' \text{ or } k_n' \leq k_n < k_n''}} +
\sum_{ \substack{k_2, \ldots, k_{n-2} \geq 0 \\ k_{n-1} \geq k_{n-1}' \text{ and } k_n \geq k_n''}} \right)
\frac{q^{u(k_2, \ldots, k_n)} T^{d_2 k_2 + \cdots + d_n k_n}}{1 - q^{g_1(k_2
+ \cdots + k_n, \ldots, k_n)} T^{d_1}}
\]
Notice that the first summation takes on a similar form to
$F_1(T)$ but with $n-3$ variables. Hence, it is meromorphic by
the induction hypothesis. Let us consider the second summation
\[
F_3(T) := \sum_{ \substack{k_2, \ldots, k_{n-2} \geq 0 \\ k_{n-1} \geq k_{n-1}' \text{ and } k_n \geq k_n''}}
\frac{q^{u(k_2, \ldots, k_n)} T^{d_2 k_2 + \cdots + d_n k_n}}{1 - q^{g_1(k_2
+ \cdots + k_n, \ldots, k_n)} T^{d_1}}.
\]
where the degree of $k_{n-2}$ in $u$ is at least one and the
degree of $k_{n-1}$ and $k_n$ in $u$ is at least two, and the
coefficients of their largest powers always take positive
values. Suppose the degree of $k_{n-2}$ is one. Then $g_1$ must
be a nonnegative constant and $g_0$ cannot contain any $k_2,
\ldots, k_{n-3}$. Thus $F_3(T)$ takes the form
\[
\sum_{ \substack{k_2, \ldots, k_{n-2} \geq 0 \\ k_{n-1} \geq k_{n-1}' \text{ and } k_n \geq k_n''}}
\frac{{ q^{c_1(k_2 + \cdots + k_n) + c_2 k_{n-2}
+ h(k_{n-1}, k_n)}} T^{d_2 k_2 + \cdots + d_n k_n}}{1 - q^{c_1} T^{d_1}}
\]
where $c_1$ and $c_2$ are constants. From geometric series and
Lemma \ref{L: MeroLemma} below, this is a meromorphic function.
Thus, we may suppose the $u$ in $F_3(T)$ has the degrees of
$k_{n-2}, k_{n-1},$ and $k_n$ being at least two, and the
coefficients of their largest powers are always positive.

Continuing this argument, we are left with a series of the form
\[
\sum_{k_2 \geq k_2', \ldots, k_n \geq k_n'}
\frac{q^{u(k_2, \ldots, k_n)} T^{d_2 k_2 + \cdots + d_n k_n}}{1 - q^{g_1(k_2
+ \cdots + k_n, \ldots, k_n)} T^{d_1}}.
\]
where $deg_{k_i}(u) \geq 2$ for each $i = 2, \ldots, n$ and the
leading coefficient of $k_i$ in $u$ is a polynomial in
$k_{i+1}, \ldots, k_n$ which takes positive values for all
$k_{i+1} \geq k_{i+1}', \ldots, k_n \geq k_n'$. It follows from
Lemma \ref{L: MeroLemma} below that this is $p$-adic
meromorphic, finishing the proof of the result.
\end{proof}

\begin{lemma}\label{L: MeroLemma}
Let $u(k_1, \ldots, k_n)$ be an increasing polynomial with
$deg_{k_i}(u) \geq 2$ for every $i = 1, \ldots, n$, and the
leading coefficient of $k_i$ in $u$ takes positive values for
all $k_1, \ldots, k_n \geq 0$. Let $v(k_1, \ldots, k_n)$ be an
increasing polynomial. Then
\begin{equation}\label{E: MeroSum}
\sum_{k_1, \ldots, k_n \geq 0} \frac{q^{u(k_1, \ldots, k_n)}
T^{d_1 k_1 + \cdots + d_n k_n}}{1 - q^{v(k_1 + \cdots + k_n, \ldots, k_n)} T}
\end{equation}
is $p$-adic meromorphic.
\end{lemma}

\begin{proof}
We will proceed by writing (\ref{E: MeroSum}) as a quotient of
two entire functions. First, the denominator of (\ref{E:
MeroSum}) takes the form
\[
G(T) := \prod (1-q^{v(k_1 + \cdots + k_n, \ldots, k_n)} T)
\]
where the product runs over all $k_1, \ldots, k_n \geq 0$ with
distinct values under $v$; that is,
\[
v(k_1' + \cdots + k_n', \ldots, k_n') \not= v(k_1 + \cdots + k_n, \ldots, k_n)
\]
whenever $k_i' \not= k_i$ for any $i$. This will be an entire
function by Lemma \ref{L: entireDenom} below, and its Newton
polygon will be bounded below by a quadratic polynomial. That
is, if we write
\[
G(T) = \sum_{r\geq 0} M_r T^r
\]
then there exists $c > 0$ and $d$ such that
\[
ord_p(M_r) \geq c r^2 + d \quad \text{ for all } r \geq 0.
\]

Let us now consider the numerator of (\ref{E: MeroSum}), which
takes the form
\begin{equation}\label{E: numerator}
\sum_{k_1, \ldots, k_n \geq 0} q^{u(k_1, \ldots, k_n)} T^{d_1
k_1 + \cdots + d_n k_n} G_{k_1, \ldots, k_n}(T)
\end{equation}
where
\[
G_{k_1, \ldots, k_n}(T) := G(T) / (1 - q^{v(k_1 + \cdots + k_n, \ldots, k_n)}).
\]
Since $G_{k_1, \ldots, k_n}$ is a factor of $G(T)$, its Newton
polygon is bounded below by the same quadratic polynomial as
that of $G(T)$.

Write
\[
u(k_1, \ldots, k_n) = g_r(k_2, \ldots, k_n) k_1^r + \cdots + g_0(k_2, \ldots, k_n).
\]
By hypothesis, $r \geq 2$ and $g_r$ takes positive values for
all $k_2, \ldots, k_n \geq 0$. Let $a_1$ satisfy $0 < a_1 <
g_r(0,\ldots, 0)$. Then there exists a constant $b_1$ such that
\[
u(k_1, \ldots, k_n) \geq a_1 k_1^2 + b_1 \quad \text{ for all } k_1, \ldots, k_n \geq 0.
\]
Find similar constants $a_i$ and $b_i$ for each variable $x_i$.
Let $a > 0$ be the minimum of the $a_i$ and $b$ be the minimum
of the $b_i$. Then, for each $i$
\[
u(k_1, \ldots, k_n) \geq a k_i^2 + b \quad \text{ for all } k_1, \ldots, k_n \geq 0.
\]

Now, for each $G_{k_1, \ldots, k_n}$ write
\[
G_{k_1, \ldots, k_n}(T) = \sum_{r \geq 0} M_r^{(k_1, \ldots, k_n)} T^r.
\]
As mentioned above, the Newton polygon is uniformly bounded
below:
\[
ord_p( M_r^{(k_1, \ldots, k_n)} ) \geq c r^2 + d \quad \text{ for all } r \geq 0 \text{ and } k_1, \ldots, k_n \geq 0.
\]
Write the numerator (\ref{E: numerator}) as a series $\sum_{s
\geq 0} A_s T^s$ where
\[
A_s = \sum_{r=0}^s \sum q^{u(k_1, \ldots, k_n)} M_{s-r}^{(k_1, \ldots, k_n)}
\]
the second sum running over all $d_1 k_1 + \cdots d_n k_n = s$.
Now, if $d_1 k_1 + \cdots d_n k_n = s$ then there exists $k_i$
such that $k_i \geq \frac{s}{d_i n} \geq \frac{s}{d' n}$ where
$d' := \max\{d_1, \ldots, d_n\}$. It follows that
\[
ord_p(A_r) \geq \min_{0 \leq r \leq s}\left\{ a \left( \frac{s}{d'n} \right)^2 + b + c(s-r)^2  + d \right\}
\geq \epsilon_1 r^2 + \epsilon_2
\]
for some $\epsilon_1 > 0$ and $\epsilon_2 \in \bb R$. This
shows the numerator of (\ref{E: MeroSum}) is entire, finishing
the proof.
\end{proof}

\begin{lemma}\label{L: entireDenom}
Let $v(x_1, \ldots, x_n)$ be an increasing polynomial with
$deg_{x_i}(v) \geq 1$ for each $i$. Consider
\[
G(T) := \prod (1-q^{v(k_1 + \cdots + k_n, \ldots, k_n)} T)
\]
where the product runs over all $k_1, \ldots, k_n \geq 0$ such
that
\[
v(k_1' + \cdots + k_n', \ldots, k_n') \not= v(k_1 + \cdots + k_n, \ldots, k_n)
\]
whenever $k_i' \not= k_i$ for any $i$. Writing $G(T) = \sum_{r
\geq 0} M_r T^r$ there exists constants $c > 0$ and $d$ such
that
\begin{equation}\label{E: NPlemma}
ord_p(M_r) \geq c r^2 + d \quad \text{ for all } r \geq 0.
\end{equation}
\end{lemma}

\begin{proof}
We will proceed by induction on the number of variables $n$.
Suppose $n = 1$. Write $v(x) = a_m x^m + \cdots + a_0$ where
$a_m > 0$ and $m \geq 1$. Let $c$ satisfy $0 < c < a_m$. Find
$x' \geq 0$ such that for all $x \geq x'$, we have $v(x) - c x
\geq 0$. Let $d := \min\{ v(x) - c x : 0 \leq x < x' \}$. Then
$v(x) \geq cx + d$ for all $x \geq 0$. Replacing $x$ with $k$,
we see that
\[
M_r = \sum (-1)^r q^{v(k^{(1)}) + \cdots + v(k^{(r)})}
\]
where the sum runs over all nonnegative integers $k^{(1)},
\ldots, k^{(r)}$ such that their values under $v$ are distinct.
Using that $v$ is an increasing polynomial and (\ref{E:
NPlemma})
\[
ord_p(M_r) \geq v(0) + v(1) + \cdots + v(r-1) \geq \frac{c r(r-1)}{2} + dr
\]
which proves the result for $n=1$.

Let us assume the result whenever there are fewer than $n$
variables. By Lemma \ref{L: technicalLemma}, $deg_{k_i}(v) \geq
1$ for every $i$. As we did for $u$ in the proof of Lemma
\ref{L: MeroLemma}, find constants $c > 0$, $d \in \bb R$, and
$k_1', \ldots, k_n' \geq 0$ such that
\begin{equation}\label{E: lowerbound}
v(k_1 + \cdots + k_n, \ldots, k_n) \geq c k_i + d \quad \text{ whenever $k_i \geq k_i'$ for all $i$.}
\end{equation}
Write
\[
G(T) = G_1(T) G_2(T)
\]
where
\[
G_1(T) = \prod_{\substack{ k_1, \ldots, k_n \geq 0, \text{ and} \\ \exists i \text{ such that } 0 \leq k_i < k_i' \\
\text{(distinct $v$ values)}}}
(1- q^{v(k_1 + \cdots + k_n, \ldots, k_n)} T)
\]
and
\[
G_2(T) = \prod_{ \substack{k_i \geq k_i' \text{ for every } i = 1, \ldots, n \\ \text{(distinct $v$ values)}}}
(1- q^{v(k_1 + \cdots + k_n, \ldots, k_n)} T).
\]
By the induction hypothesis, $G_1(T)$ is an entire function
whose Newton polygon is bounded below by a quadratic
polynomial. For convenience, let us make the change of
variables $k_i \mapsto k_i - k_i'$ in $G_2(T)$. Thus, consider
\[
G_2(T) = \prod_{ \substack{k_1, \ldots, k_n \geq 0 \\ \text{(distinct $v$ values)}}}
(1- q^{v(k_1 + \cdots + k_n, \ldots, k_n)} T).
\]
Notice that the set
\[
\cc A := \{ (k_1 + \cdots + k_n, \ldots, k_n) : k_i \geq 0\} = \{ (a_1, \ldots, a_n) : a_1 \geq \cdots \geq a_n \geq 0 \}
\]
is linearly ordered by lexicographic ordering, and has a
minimal element $(0, \ldots, 0)$. For example, when $n=3$ we
have
\[
(0,0,0) < (1,0,0) < (1,1,0) < (1,1,1) < (2,0,0) < \cdots.
\]
It follows that if $a^{(0)} < a^{(1)} < \cdots < a^{(r)}$,
where
\[
a^{(i)} = (a_1^{(i)}, \ldots, a_n^{(i)}) \in \cc A,
\]
then $a_1^{(i)} \geq \frac{i}{n}$ for each $i$. Now, write
$G_2(T) = \sum_{r \geq 0} \overline{M}_r T^r$ with
\[
\overline{M}_r := \sum (-1)^r q^{v(a^{(1)}) + \cdots v(a^{(r)})}
\]
where the sum runs over all $a^{(i)} \in \cc A$ satisfying
$v(a^{(i)}) \not= v(a^{(j)})$ for $i \not= j$. Ordering the
$a^{(i)}$ if necessary we may assume $a_1^{(i)} \geq
\frac{i-1}{n}$ for every $i$. Since
\[
a_1^{(i)} = k_1^{(i)} + \cdots + k_n^{(i)},
\]
we see that there must be a $j$ such that $k_j^{(i)} \geq
\frac{i-1}{n^2}$. Consequently, from (\ref{E: lowerbound}) we
see that
\[
ord_p(\overline{M}_r) \geq \sum_{i=1}^r \left( \frac{c(i-1)}{n^2} + d \right) = \frac{c r (r-1)}{2 n^2} + dr
\quad \text{for all } r \geq 0.
\]
This proves the result for $G_2(T)$. Since the result is true
for both $G_1$ and $G_2$, it is true for their product.
\end{proof}

This finish the proof of Theorem \ref{T: MeroResult}.

\bibliographystyle{amsplain}
\bibliography{XBib}

\end{document}